\documentclass[fleqn]{mat01}
\usepackage{times,mathtimy,amssymb,latexsym}
\begin{document}

\setcounter{page}{389}
\firstpage{389}

\newtheorem{rema}{Remark}
\newtheorem{lemm}{Lemma}
\newtheorem{theore}{\bf Theorem}

\def\theor{\trivlist\item[\hskip\labelsep{\bf Theorem.}]}
\def\lem{\trivlist\item[\hskip\labelsep{\it Lemma.}]}
\def\defini{\trivlist\item[\hskip\labelsep{\rm DEFINITION.}]}
\def\pot{\trivlist\item[\hskip\labelsep{\it Proof of Theorem $1$.}]}

\title{Some remarks on good sets}

\markboth{K~Gowri Navada}{Some remarks on good sets}

\author{K~GOWRI NAVADA}

\address{Department of Mathematics, University College, Mangalore
University, Hampanakatta~575~001, India\\
\noindent E-mail: gnavada@yahoo.com}

\volume{114}

\mon{November}

\parts{4}

\Date{MS received 3 February 2004; revised 4 July 2004}

\begin{abstract}
It is shown that (1) if a good set has finitely many related components,
then they are full, (2) loops correspond one-to-one
to extreme points of a convex set. Some other properties of good sets
are discussed.
\end{abstract}


\maketitle

\section*{Introduction and preliminaries}

In this note we make some remarks on good sets in $n$-fold
Cartesian product as defined in \cite{KNB}. We need the following
definitions:\vspace{.5pc}

Let $X_1,X_2,\ldots, X_n$ be non-empty sets and let $\Omega =X_1\times
X_2\times \cdots\times X_n$ be their Cartesian product. For each
$i,1\leq i \leq n$, $\Pi_i$ will denote the canonical projection of
$\Omega$ onto $X_i$. A subset $S \subset \Omega$ is said to be {\it
good}, if every complex valued function $f$ on $S$ is of the form:
\begin{equation*}
f(x_1,x_2,\ldots,x_n) = u_1(x_1) + u_2(x_2) + \cdots + u_n(x_n),
(x_1,x_2,\ldots,x_n) \in S,
\end{equation*}

$\left.\right.$\vspace{-1.5pc}

\noindent for suitable functions $u_1,u_2,\ldots, u_n$ on
$X_1,X_2,\ldots,X_n$ respectively (\cite{KNB2}, p.~181).

A subset $S$ of $\Omega$ is said to be {\it full}, if $S$ is a maximal
good set in $\Pi_1 S\times \Pi_2 S\times\cdots\times \Pi_n S$
(\cite{KNB2}, p.~183).

Two points $ x,\ y$ in a good set $S$ are said to be {\it related},
denoted by $xRy$, if there exists a finite subset of $S$, which is full
and contains both $x$ and $y$. $R$ is an equivalence relation, whose
equivalence classes we call as the {\it related components} of $S$. Note
that related components of $S$ are full subsets of $S$ (\cite{KNB2}, p.~190). 

\begin{rema}
{\rm Here we prove:}
\end{rema}

\begin{theore}[\!]
If a full set $F$ has finitely many related components{\rm :}
$F=\cup_{i=1}^k R_i${\rm ,} then $k=1$.\vspace{.7pc}
\end{theore}

To  prove this we need some preliminary results.

Let $S$ be a good set, $S=\cup R_\alpha$ where $R_\alpha$ are its
related components. Define an equivalence relation $E_i$ on
$\Pi_i S$ as follows: $x_{i}E_iy_{i}$ if there exists a finite sequence
$R_1, R_2, \ldots, R_k$ such that $x_{i} \in \Pi_i R_1, y_{i} \in \Pi_i
R_k$ and $\Pi_i R_j\cap \Pi_i R_{j+1} \neq \emptyset$ for $1 \leq j \leq
k- 1 $ (\cite{KNB2}, p.~189). For $x_{i}\in \Pi_i S$, $[x_{i}]$ denotes
the $E_i$-equivalence class of $x_{i}$. If an element in $\Pi_iR_\alpha$
is $E_i$-equivalent to an element in $\Pi_i R_\beta$, we will say that
$R_\alpha$ and $R_\beta$ are \hbox{$E_i$-equivalent.} Let $C$ be a
cross-section of $R_\alpha$'s. Let ${\cal F}_i$ be the set of all
$E_i$-equivalence classes for $1\leq i\leq n$. Define  $\phi \hbox{:}\
C \rightarrow { \cal F}_1 \times {\cal F}_2 \times \cdots \times {\cal
F}_n$ by $f(x_1,\ldots, x_n)=([x_1],\ldots, [x_n])$. 

\begin{lemm}{\it 
$\phi(C)$ is good. $S$ is full if and only if $\phi(C)$ is full.}
\end{lemm}

\begin{proof}
The $\phi\hbox{:}\ C \mapsto \phi(C)$ is one-to-one: If $(x_1,\ldots, x_n)
\neq (y_1,\ldots, y_n)$ are in $C$ then they belong to two different
related components. If $([x_1],\ldots ,[x_n])=([y_1],\ldots ,[y_n])$
then these two related components, say $R_{1}$ and $R_{2}$, have all the
coordinates equivalent which is not possible. To prove this, let us
define a function $h$ which is equal to zero everywhere in $S$ except on
$R_1$ where it is a non-zero constant. There exists functions $u_i$
defined on $\Pi_i S$, $1 \leq i \leq n$, such that $u_1+\cdots +u_n = h$
on $S$. As the function $h$ is constant on each $R_\alpha$, the $u_i$
are constants on $\Pi_i R_\alpha$ (\cite{KNB2}, p.~185, Corollary~2).
Let $ c_1,\ldots ,c_n$ be these constants $\Pi_{i}R_{2}, 1\leq i\leq n$.
As all the coordinates of $ R_1$ and $ R_2$ are equivalent we get the
same constants $c_1,\ldots c_n$ on the coordinates of $R_1$. But
$c_1+\cdots +c_n=0$ and $h$ is non-zero on $R_1$. This contradicts the
fact that $h = u_1 + u_2 +\cdots +u_n$. Therefore $\phi$ is one-to-one. 

Next we show that $\phi(C)$ is good. Take a function $h$ on $\phi(C)$. This
defines in a natural manner a function on $C$. Denote it also by $h$.
Define $g$ on $S$ by taking it as a constant on each $R_\alpha$, i.e.,
$g(y_1,\ldots, y_n) = h(x_1,\ldots,x_n)$ for all $(y_1,\ldots ,y_n)\in
R_\alpha$ where $(x_1,\ldots,x_n)$ is in $C\cap R_\alpha $. There exists
$u_1,\ldots, u_n$ such that $u_1+\cdots +u_n = g$ on $S$ as $S$ is good.
Since $g$ is constant on each $ R_\alpha$, $ u_i$ is constant on $\Pi_i
R_\alpha$. Define $ v_i([x_i])=u_i(x_i)$ for all $x_i, 1\leq i\leq n$.
These functions are well-defined. Further, 
\begin{align*}
v_1([x_1])+\cdots +v_n([x_n]) &= u_1(x_1)+\cdots +u_n(x_n)\\[.2pc]
&= g(x_1, \ldots,x_n) = h([x_1],\ldots, [x_n]).
\end{align*}
This shows $\phi(C)$ is good.

Suppose $S$ is full. If $\phi(C)$ is not full,
then given the zero function on $\phi(C)$ there exist two distinct sets of
functions $\{v_i\}$ and $\{v_i'\}$ defined on the $i$th coordinate space
of $\phi(C)$ for each $i$ whose sum is equal to zero which also satisfy
\begin{equation*}
v_i([x_i^0])=v_i^{ \prime}([x_i^0])=0
\end{equation*}
for some $[x_i^0]\in\Pi_i(\phi(C))$ for $1\leq i\leq n-1$. Define $\{u_i\}$
and $\{u_i'\}$ on the $i$th coordinate space $\Pi_iS$ by 
\begin{equation*}
u_i(x_i)=v_i([x_i])\ \hbox{and} \ u_i'(x_i)=v_i'([x_i])
\end{equation*}
for $1\leq i\leq n$. Then the sum of $u_i$ as
well as $u_i'$ is equal to zero but they are different solutions (with
the same boundary conditions) because $v_i$ and $v_i'$ are different.
This contradicts the fact that $S$ is full.

Conversely, if $\phi(C)$ is full, then we prove $S$ is also full. For this, take the zero
function on $S$. Suppose there are two distinct sets of functions
$\{u_i\}$ and $\{u_i'\}$ with 
\begin{equation*}
u_1+\cdots +u_n = u_1'+\cdots +u_n' =0
\end{equation*}
on $S$ with $u_i(x_i^0)=u_i'(x_i^0)=0$ for some $x_i^0\in \Pi_i(S),\
1\leq i\leq n-1$. All the functions $u_i$ and $u_i'$ are constant on
each $\Pi_i R_\alpha$ so also on each $E_i$ equivalence class. Define
$v_i([x_i])=u_i(x_i)$ and $v_i'([x_i])=u_i'(x_i), \forall i$. Then
$\{v_i\}$ and $\{v_i'\}$ are distinct solutions of the zero function on
$\phi(C)$ which also satisfy $v_i([x_i^0])=v_i'([x_i^0])=0$ for $1\leq
i\leq n-1$. Since $\phi(C)$ is full, the functions $v_i$ and $v_i'$ are the
same which implies $u_i$ and $u_i'$ are equal. So $S$ is full.\hfill
$\Box$
\end{proof}

\begin{defini}$\left.\right.$\vspace{.5pc}

\noindent {\rm For a finite good set $S$, we call the cardinality of
$\cup_{i=1}^n\Pi_iS$ as the {\it number of coordinates} of $S$ and denote it
by $N(S)$. The cardinality of $\Pi_iS$ is called the {\it number of} $i$-{\it th
coordinates}\break of $S$.}
\end{defini}

A finite good set $S$ is full if and only if $N(S)-(n-1)=|S|$. If $S$ is
good then $N(S)-(n-1)\geq |S|$. For a finite set $S$ if $|S|>N(S)-(n-1)$ then
$S$ is not good (\cite{KNB}, p.~80).

\begin{pot}
Suppose $F= \cup_{i=1}^kR_i$ is full. We want to show that there is a
finite, full subset $S$ of $ F$ which intersects each $R_i$. If $k>1$
this will be a contradiction to the fact that $R_i$'s are related
components. Since $F$ is full, by lemma~1, $\phi(C)$ is full. It has $k$
points and dimension $n$. So the total number of coordinates in $\phi(C)$
is $k+(n-1)$. Let these coordinates be labeled as
$\alpha_1,\ldots,\alpha_{k+n-1}$ in some order. There are $k$ points
each having $n$ entries and each of these $nk$ entries should be one of
these $k+n-1$ coordinates of $\phi(C)$. So we get a partition of $nk$ as
$nk=l_1+l_2+\cdots +l_{(k+(n-1))}$ where $l_i$ denotes the number of
times $\alpha_i$ is repeated. When a coordinate, say $[x_j] \in \Pi_j
\phi(C)$, is repeated in $\phi(C)$, it means the corresponding two related
components of $F$ are $E_j$-equivalent. If $[x_j]$ occurs $l$ times in
$\phi(C)$ then $l$ number of related components of $F$ are $E_j$-equivalent. For this, it is necessary to have at least $l-1$ different
pairs of related components $(R_\alpha, R_\beta)$ such that
$\Pi_jR_\alpha\cap \Pi_jR_\beta\neq \emptyset$: 

Suppose $l$ related components are $E_j$-equivalent. Consider a graph
whose vertices are these related components and whose edges are pairs of
related components among these which have at least one common $j$th
coordinate. This graph is connected because the related components are
$E_j$-equivalent. The number of vertices is $l$ so there should be at
least $l-1$ edges in it.

In this way we get totally (at least) $l_1-1+l_2-1+\cdots
+l_{(k+(n-1))}-1=nk-(k+(n-1))$ pairs of related components $(R_\alpha,
R_\beta)$ such that for some $i,$ $\Pi_iR_\alpha\cap
\Pi_iR_\beta\neq\emptyset$. For each such pair $(R_\alpha, R_\beta)$
take one point from each of the two related components $R_\alpha $ and
$R_\beta$ such that the chosen points have the same $i$th coordinate.
All these points together form a finite subset of $F$. The intersection
of this set with each $R_i$ is also finite and non-empty. (Note that
since $F$ is full, each $R_\alpha$ has a common coordinate with some
other related component.) Take the finite full set $F_i\subset R_i$
which contains this intersection. (Any finite subset of a related
component is contained in a finite full set.) Let $ S=\cup_{i=1}^k F_i$.
Then $S$ is a finite subset of $F$.

To show that $S$ is full we have to find the number of coordinates of
$S$ and the number of points in $S$. Let $A_i$ denote the number of
coordinates of $F_i$. Then, since $F_i$ is full, the number of points in
$F_i$ is $A_i-(n-1)$. So the number of points in $S$ is $|S|=A-k(n-1)$
where $A=A_1+\cdots +A_k$. Now the number of coordinates of $S$ is no
more than $A$. In this counting, if $F_\alpha$ and $F_\beta$ have a
common coordinate, then this common coordinate will be counted once each
in $A_\alpha$ and $A_\beta$. But we know that there are at least
$nk-(k+(n-1))$ such pairs $F_\alpha, F_\beta$. So the number of
coordinates of $S$ is at most
\begin{equation*}
A-(nk-(k+(n-1)))= A-(n-1)(k-1)=|S|+(n-1).
\end{equation*}
But the number of coordinates of $S$ cannot be lesser than
this: if it is the case $S$ will not be good. This shows the number
of coordinates of $S$ is equal to $|S|+(n-1)$. So $S$ is full. If $k>1$,
this is a contradiction as noted at the beginning of the proof.\hfill $\Box$
\end{pot}

\begin{rema}
{\rm Here we show the connection between loops and extreme points of a
convex set. We need the following definitions.}
\end{rema}

\pagebreak

Given any finitely many symbols $t_1,t_2,\ldots, t_k$ with {\it
repetitions} allowed and given any finitely many integers
$n_1,n_2,\ldots, n_k$, we say that the formal sum $n_1t_1 + n_2t_2 +
\cdots + n_kt_k$ {\it vanishes}, if for every $t_j$ the sum of the
coefficients of $t_j$ is equal to zero (\cite{KNB2},\break p.~183).

\begin{defini}$\left.\right.$\vspace{.5pc}

\noindent Let $\Omega = X_1\times X_2\times\cdots\times X_n$. A
non-empty finite subset $ L = \{x_1,x_2, \ldots,x_k\}$ of $\Omega $ is
called a {\it loop}, if there exist non-zero integers
$n_1,n_2,\ldots,n_k$ such that the sum $\sum_{j=1}^kn_jx_j$ vanishes in
the sense that the formal sum vanishes coordinate-wise and no strictly
smaller non-empty subset of $L$ has this property (\cite{KNB2}, p.~183).
\end{defini}

\begin{lemm}
Let $L=\{x_1,\ldots ,x_k\}$ be a loop. Then there is a unique
{\rm (}except for the sign{\rm )} set of integers $n_1,\ldots, n_k$ with
$gcd(n_1, \ldots,n_k)=1$ such that the formal sum\break $\sum_{j=1}^kn_jx_j$ vanishes.
\end{lemm}

\begin{proof}
Suppose there are two sets of integers $\{n_j\}$ and $\{m_j\}$ with
these properties. Also assume there is a $p$ for which $|n_p|\neq
|m_p|$. Then $\sum_{j=1}^kn_jx_j=0$ and $\sum_{j=1}^km_jx_j=0$ imply
$\sum_{j=1 }^k(m_pn_j-n_pm_j)x_j=0$, where the co-efficient of $x_p$
vanishes. As $L$ is a loop any proper subset of $L$ is not a loop. So
$m_pn_j=n_pm_j$ for $1\leq j\leq n$. Then $gcd(m_pn_1,\ldots, m_pn_k)
= gcd(n_pm_1,\ldots ,n_pm_k)$, i.e., $|m_p|=|n_p|$. We have to
prove that either for all $j,$ $m_j=n_j$ or for all $j, m_j=-n_j$.
Suppose $m_j=n_j$ for some $l$ number of $j$'s for $0<l<n$ and
$m_j=-n_j$ for the remaining $n-l$ number of $j$'s. Then adding the
equations $\sum_{j=1}^kn_jx_j=0$ and $\sum_{j=1}^km_jx_j=0$ we get a
smaller formal sum (containing only $l$ terms) to be zero which is a
contradiction to the minimality of the loop $\sum_{j=1}^kn_jx_j=0$.
\hfill $\Box$
\end{proof}

\begin{lemm}
Let $L=\{x_1,\ldots ,x_k\}$ be a loop. Let $n_j$ be as in Lemma~{\rm 2} above.
Suppose the formal sum $\sum_{j=1}^kr_jx_j=0$ for some real numbers
$r_j$. Then there exists a real number $\alpha$ such that $r_j=\alpha
n_j$ for each $j$.
\end{lemm}

\begin{proof}
If we assume that $r_j$'s are rationals, then the result is easy to prove.

To prove the general case, note that the formal sum
$\sum_{j=1}^kr_jx_j=0$ gives a set of $N$ homogeneous equations -- one for
each coordinate in $\cup_{i=1}^n\Pi_iL$ where \hbox{$N=|\!\!\cup_{i=1}^n\Pi_iL|$.}
The matrix corresponding to this set of equations gives a linear map
from $R^k$ to $R^N$. This matrix consists only of 0's and 1's and so can be
thought of as a linear map from $Q^k$ to $Q^N$ the kernel of which is
one-dimensional (by the result for the rational case). It means that the
rank of the matrix is $k-1$. This is also the rank of the matrix, when
the linear map is considered from $R^k$ to $R^N$. The null space of this
matrix is one-dimensional, i.e., there exists some $\alpha$ such that
$r_j=\alpha n_j$ for $1\leq j\leq k$.\hfill $\Box$
\end{proof}

This also shows that $\sum_{j=1}^mr_jx_j=0$ is not possible for $m<k$
even for real $r_j$'s. The above proof in fact shows that a finite set
$\{x_1,x_2,\ldots, x_k\}$ of points in $X_1\times X_2\times \cdots
\times X_n$ is a loop if and only if there is a finite set of non-zero
real numbers $r_1, r_2, \ldots, r_k$ such that the formal sum
$\sum_{j=1}^kr_jx_j$ vanishes and no proper subset of $\{x_1, x_2,
\ldots, x_k\}$ has this property. If we assume that $\sum_{j=1}^k |
r_j|=1$ then $|\alpha|=(\sum_{j=1}^k|n_j|)^{-1}$.

Let $S\subset X_1\times\cdots\times X_n$ be a finite set, not
necessarily good. Let $C(S)$ be the set of all functions on $S$.
The norm in $C(S)$ is defined by $\|f\|= \max_{x\in S}|f(x)|$. 
Let $U(S)=\{f\in C(S) |f(x_1,\ldots, x_n)=u_1(x_1)+\cdots +u_n(x_n)$ where
$u_i$ is a function on $\Pi_i(S)\}$. $U(S)$ is a subspace of $C(S)$. Let
$M(S)$ denote the space of all signed measures on $S$ with the total
variation norm (which is just the $L_1$ norm). 
Then $M(S)=(C(S))^*$. Take the
subspace $(U(S))^\perp \subset M(S)$. This is the set of all signed measures $\mu$ with
$\mu(f)=0, \forall f\in U(S)$. Note that $\mu \in (U(S))^\perp$ if and
only if all the one-dimensional marginals of $\mu$ vanish. Consider
\begin{equation*}
A=\{\mu\in (U(S))^\perp:\|\mu\|\leq 1\}. 
\end{equation*}
This set is convex.

\begin{defini}$\left.\right.$\vspace{.5pc}

\noindent By a {\it weak loop} we mean a finite set $\{x_1, x_2, \ldots, x_l\}\in
X_1\times X_2 \times \cdots \times X_n$ for which there exist real numbers
$r_1, r_2, \ldots, r_l$, with at least one $r_i$ non-zero, such that the
formal sum $\sum_{i=1}^lr_ix_i = 0$, coordinate-wise.
\end{defini}

\begin{theore}[\!]
The extreme points of $A$ are given by $\mu_L$ where $L=\{x_1,\ldots,
x_k\}$ is a loop. In this case
\begin{equation*}
\mu_L(x_j)=n_j(|n_1|+\cdots  +|n_k|)^{-1},
\end{equation*}
where $(n_1,\ldots, n_k)$ are given by $\sum_{j=1}^k n_jx_j=0${\rm ,}
and for all other $x\in S, \mu_L(x)=0$.
\end{theore}

\begin{proof}
First we note that $\mu_L\in (U(S))^\perp\hbox{:}$ It is enough to show that
$\mu_L(u_i)=0$ where $u_i$ is a function on $\Pi_iS$, and this is easily
verified from the form of $\mu_L$ and the fact that $L$ is a loop. To
show that $\mu_L$ is an extreme point of $A$ suppose
$\mu_L=a\lambda+b\nu$ where $a+b=1, a>0,b>0$ and $\lambda,\nu \in A$.
Then $\|\lambda\|, \|\nu\| = 1$ and $\lambda,
\nu \in (U(S))^\perp$. Restricting all these three measures to $L$, we have
$\mu_L=\mu_L|_L=a\lambda|_L+b\nu|_L$, and since $\|\mu_L\|
=1$, we have $\|\lambda|_L\| = \|\nu|_L\| =1$.
This shows that $\lambda$ and $\nu$ are supported on $L$. Denote
$\lambda(x_j)=\alpha_j, \nu(x_j)=\beta_j, 1\leq j\leq k$. Since
$\lambda$, $\nu$ are in $U(S)^\perp$, their marginals vanish, which is
equivalent to saying that the formal sum $\sum_{j=1}^k\alpha_jx_j=0,
\sum_{i=1}^k\beta_jx_j =0$. Since $\sum_{j=1}^k |\alpha_j| = 1,\
\sum_{j=1}^k|\beta_j| =1$. By the above lemma $\lambda = + \nu$ or
$\lambda = -\nu$, and since $a, b >0$ we see that $\mu_L = \lambda =
\nu$, and $\mu_L$ is an extreme point of $A$.

Conversely, take an extreme point $\mu$ of $A$. Then $\|\mu\|=1$. We
show that support of $\mu$ is a weak loop. Let $\{x_1,\ldots ,x_k\}$ be
the support of $\mu$. Since $\mu\in (U(S))^\perp$, for any $i$, if $u_i$
is a function on $\Pi_i(S)$ then $\mu(u_i)=0$, i.e., $\sum_{j=1}^k\mu
(x_j)u_i(x_j)=0$. Take $u_i(x_{1i})=1$ and $u_i(x)=0$ for all {\rm
other} ~$x\in \Pi_i(S)$. Then $\mu(u_i)=\sum\mu(x_j)=0$ where the sum
runs over all $x_j$ for which $x_{ji}=x_{1i}$. With similar arguments
for other $x_{ji}$ we see that the formal sum
$\sum_{j=1}^k\mu(x_j)x_j=0$. 

Now we prove that this weak loop has to be a loop. Call $\mu(x_j)$ as
$m_j$. Then we have $\sum_{j=1}^km_jx_j=0$. Suppose this is not a loop.
We prove that the measure $\mu$ is not an extreme point of $A$. Any weak
loop contains a loop. Let $\sum n_jx_j=0$ be this loop. Here the sum
runs over a proper subset of $\{1,\ldots ,k\}$. Then taking $n_j=0$
whenever necessary $\sum_{j=1}^km_jx_j=\sum_{j=1}^k n_jx_j+\sum_{j=1}^k
(m_j-n_j)x_j$ is a sum of two weak loops. Note that the two weak loops
on the right side are not multiples of each other. Let $\mu_1, \mu_2$ be
the measures corresponding to $\sum_{j=1}^k n_jx_j$ and $\sum_{j=1}^k
(m_j-n_j)x_j$ respectively. That is, $\mu_1(x_j)=n_j,\
\mu_2(x_j)=m_j-n_j$ for $1\leq j\leq k$ and for any other $x,\ \mu_1(x)=
\mu_2(x)=0$. Clearly $\mu_1, \mu_2\in (U(S))^\perp $. Then
$\mu_1+\mu_2=\mu$ and $\|\mu\|=1\leq \|\mu_1\|+\|\mu_2\|$. If for each
$j$ the coefficients $n_j$ and $m_j-n_j$ of $x_j$ on the right-hand side
have the same sign, then $|m_j|=|n_j|+|m_j-n_j|$ for each $j$ so that
$\|\mu\|=1= \|\mu_1\|+\|\mu_2\|$. Then we can write
\begin{equation*}
\sum_{j=1}^km_jx_j=\|\mu_1\| \left(\sum_{j=1}^kn_jx_j(\|\mu_1\|)^{-
1}\!\right)\!+\!\|\mu_2\| \left(\sum_{j=1}^k(m_j\!-\!n_j)x_j(\|\mu_2\|)^{-1}\!\right)
\end{equation*}

$\left.\right.$\vspace{-1.5pc}

\noindent which shows $\mu$ is not an extreme point of $A$.

Now we show that the measure which is supported on a weak loop can be
written as a sum of two measures $\mu_1$ and $\mu_2$, both in
$U(S)^\perp$, with $\|\mu\|=\|\mu_1\|+\|\mu_2\|$. We already know that
$\sum_{j=1}^km_jx_j$ can be written as 
\begin{equation*}
\sum_{j=1}^km_jx_j=\sum_{j=1}^kn_jx_j +\sum_{j=1}^k r_jx_j,
\end{equation*}
where the two weak loops on the right-hand side are not multiples of
each other. In this representation we want $|m_j|= |n_j|+|r_j|$ for each
$j$, i.e., $n_j$ and $r_j$ should have the same sign. Suppose for some
$j_0$ this does not happen. Let us assume $n_{j_0}>0$, $r_{j_0}<0$ and
$|r_{j_0}|>|n_{j_0}|$. Then we can write 
\begin{equation*}
\sum_{j=1}^km_jx_j=\sum_{j=1}^k(n_j-(n_{j_0}/r_{j_0}) r_j)x_j+
\sum_{j=1}^k(r_j+(n_{j_0}/r_{j_0}) r_j)x_j.
\end{equation*}
The right-hand side is a sum of two weak loops; the first does not
contain the term $x_{j_0}$ and the second contains $m_{j_0}x_{j_0}$. If
for some $j$, $n_j$ and $r_j$ have the same sign, then since
$(n_{j_0}/r_{j_0})<0, -(n_{j_0}/r_{j_0})r_j$ has the same sign as $r_j$
so $n_j-(n_{j_0}/r_{j_0})r_j$ has the same sign as $n_j$. Also
$r_j+(n_{j_0}/r_{j_0}) r_j$ has the same sign as $r_j$ because
$|(n_{j_0}/r_{j_0})|<1$. This shows that the new coefficients of $x_j$
in the new weak loops have the same sign if they had the same sign in
the original weak loops. Also because the original weak loops are not
multiples of each other, these two weak loops cannot have all the
coefficients equal to 0. In this way we get another representation of
the left-hand side as a sum of weak loops with lesser number of $j$'s for which
$n_j$ and $r_j$ have opposite signs. Again the two weak loops are not
multiples of each other because $x_{j_0}$ is present in the second weak
loop but not in the first. Applying the same procedure repeatedly we get
the two weak loops with all $j$ having the $n_j$ and $r_j$ of the same
sign. This proves $\mu$ is not an extreme point of $A$.\hfill $\Box$ 
\end{proof}

\begin{rema}{\rm 
Here we discuss some properties of a maximal good set contained in a
given set.}
\end{rema}

Let $S\subset X_1\times\cdots \times X_n$, $S$ not necessarily good.
Consider the collection ${\cal {G}}$ of good subsets of $S$. This
collection is closed under arbitrary increasing unions, hence by Zorn's
lemma there exists a maximal set $M$ in ${\cal{G}}$. Note that
$\Pi_i(M)=\Pi_i(S)$ for $1\leq i\leq n $. Denote the space of all
functions on $M$ as $ C(M)$. Call a function $f$ {\it good} if for each
$i$ there is a function $u_i$ on $\Pi_i S$ such that
\begin{equation*}
f(x_1,\ldots,x_n)= u_1(x_1) + u_2(x_2) +\cdots + u_n(x_n)
\end{equation*}
for all $(x_1,\ldots,x_n)\in S$. Let $ U(S)$ denote the class of good
functions on $S$.

\begin{theore}[\!]
The map $f\mapsto f|M$ is a one-to-one linear map from $U(S)$ onto
$C(M)$.
\end{theore}

\begin{proof}
Clearly this map is linear. It is also onto. To prove this, take any
function $g$ on $M$. Then $g=u_1+\cdots +u_n$. Here $u_i$ are defined on
$\Pi_i(M)$ which is same as $\Pi_i(S)$. Then $f$ defined on $S$ by
$f=u_1+\cdots +u_n$ has the property that $f|M=g$. Because $M$ is a
maximal good set, any $x\in S\backslash M$ forms a loop with some
finitely many elements of $M$. Let this loop be $\{x, y_2,\ldots ,y_k\}$
where $y_j$ are from $M$. Then $n_1x+\sum_{j=2}^k n_jy_j=0$ for some
integers $n_j$. This loop is unique because if $\{x, z_2,\ldots ,z_l\}$
where $z_j$ are from $M$ is another loop, then $m_1x+\sum_{j=2}^l
m_jz_j=0$ for some integers $m_j$. Multiplying the first equation by
$m_1$ and the second by $n_1$ and subtracting one from the other we get
a weak loop in $M$ which is not possible because $M$ is good. 

Given the zero function on $M$, there is a unique extension of this
function to the whole of $S$ which is in $U(S)$ (namely, the zero
function): Let $0=u_1+\cdots + u_n$ on $M$. By taking $f=u_1+\cdots +
u_n$ we can extend this function to $S$. Take a point $x$ in $S\backslash M$ and
the loop $\{x,y_2,\ldots, y_k\}$ it makes with the elements of $M$. Then
the formal sum $n_1x+\sum_{j=2}^kn_jy_j=0$ which gives
$n_1x_i+\sum_{j=2}^kn_{j}y_{ji}=0$ for each $1\leq i\leq n$, where $x_i$
and $y_{ji}$ denote the $i$th coordinate of $x$ and $y_j$ respectively.
This gives $x_i = -\sum_{j=2}^k(n_{j}/n_1)y_{ji}$ formally. Then
\begin{equation*}
\sum_{i=1}^nu_i(x_i)=\sum_{i=1}^nu_i\left(-\sum_{j=2}^k(n_{j}/n_1)y_{ji}\right).
\end{equation*}
We can also write
\begin{equation*}
\sum_{i=1}^nu_i(x_i)=-\sum_{j=2}^k(n_j/n_1)\sum_{i=1}^nu_i(y_{ji})
\end{equation*}
because $\sum n_j/n_1=0$ when the sum is taken over those $j$ for which
$y_{ji}$ fixed and $\neq x_i$ and $\sum -n_j/n_1=1$ when the sum is over
those $j$ for which $y_{ji}=x_i$. But $\sum_{i=1}^nu_i(y_{ji})=0$ as
$y_j$ are in $M$. So we get $\sum_{i=1}^nu_i(x_i)=0$. This shows
$f=\sum_{i=1}^nu_i=0$ on $S$. So the zero function on $M$ has a unique
extension to a function in $U(S)$. It follows that the map $f \mapsto
f|_M$ is one-to-one from $U(S)$ to $C(M)$, and the theorem is proved.

Let $C(S)$ denote the set of all functions on $S$. Let $U(S)^\perp
=\{\mu\in C(S)^*|\mu(f)=0, \forall f\in U(S)\}$.\hfill $\Box$
\end{proof}

\begin{theore}[\!]
The dimension of $U(S)^\perp$ is $|S|-|M|$ when $S$ is finite. A
basis for $U(S)^\perp$ is given by the set of $\mu_L$ where $L$ is a
loop of the form $\{x,y_2,\ldots, y_k\}$ where $x\in S\backslash M$ and $y_2,
\ldots, y_k $ are in $M$.
\end{theore}

\begin{proof}
The dimension of $C(S)$ is $|S|$ and that of $U(S)$ is $|M|$ by the
previous theorem. The space $U(S)^\perp$ is equivalent to
$(C(S)/U(S))^*$. Therefore $\dim (U(S)^\perp) = \dim(C(S))-\dim(U(S))$,
which is equal to $|S|-|M|$. Every $x\in S\backslash M$ makes a unique loop $L_x$
with suitable elements from $M$. These loops give rise to $|S|-|M|$
measures in $U(S)^\perp$ by the theorem in Remark 2. They are linearly
independent: If $\sum c_j\mu_{L_{x_j}}=0$ then $\sum
c_j\mu_{L_{x_j}\!\!}(x)=0, \forall x\in S$. Taking $x=x_i$ we get
$\mu_{L_{x_j}}(x_i)=0, \forall j\neq i$ and $\mu_{L_{x_i}}(x_i)\neq 0$.
This gives $c_i=0$. Therefore, these $|S|-|M|$ measures are linearly
independent and so form a basis for $U(S)^\perp$. This proves the
theorem.\hfill $\Box$
\end{proof}

\begin{defini}$\left.\right.$\vspace{.5pc}

\noindent A set $S\subset X_1\times\cdots\times X_n$ is called
relatively full if there exist $x_i^0 \in \Pi_iS, 1\leq i \leq n-1$ such
that any $f\in U(S) $ has a unique representation as $f=u_1+\cdots +u_n$
when we fix the value of $u_i(x_i^0), 1\leq i\leq n-1$.
\end{defini}

It is easy to see that if $S$ is relatively full then for any choice of
elements $x_i^0\in \Pi_iS, 1\leq i \leq n-1 $ the solution of $f = u_1 +
u_2 + \cdots + u_n$, $f\in U(S)$ is unique if we fix the values of
$u_i(x_i^0), 1\leq i \leq n-1$. Moreover, if the solution is unique, with
the prescribed constraints, for the zero function, then it is unique for
all functions in $U(S)$.

\begin{theore}[\!]
If $S$ is relatively full then any maximal good set $M\subset S$ is full.
\end{theore}

\begin{proof}
Take the zero function on $M$. Fix $u_i(x_i^0)=0 ~{\rm for } ~1\leq
i\leq n-1$. Let $0=\sum_{i=1}^nu_i$ on $M$. It can be uniquely extended
to $S$ as a function $f$ in $U(S)$. Then $f=\sum_{i=1}^nu_i$ on $S$ with
$u_i(x_i^0)=0$ for $1\leq i\leq n-1$. These $u_i$'s are unique
because $S$ is relatively full. Therefore $M$ is full. This proves the
theorem.

Any set $S\subset X_1\times\cdots\times X_n$ can be written uniquely as
a disjoint union $S=\cup R_\alpha$ of maximal relatively full sets
$R_\alpha$ of $S$: A one point set of $S$ is relatively full. Union of a
chain of relatively full sets is again relatively full. Using Zorn's
lemma, there exist maximal relatively full sets in $S$. As union of two
relatively full sets with non-empty intersection is again relatively
full, these maximal relatively full sets of $S$ do not intersect each
other and their union is $S$.

\begin{enumerate}
\renewcommand\labelenumi{\arabic{enumi}.}
\leftskip -.2pc
\item $R_\alpha$ and $R_\beta$ for $\alpha\neq \beta$ cannot have $n-1$
coordinates in common.

\item For $n=2 $, $\Pi_iR_\alpha\cap \Pi_iR_\beta = \emptyset$ if $\alpha\neq
\beta$.

\item In each $R_\alpha$, there exists a maximal good set $F_\alpha$
which is also full. 

\item Although each $F_\alpha$ is good and full, $\cup F_\alpha\subset
S$ need not be good, except when $n=2$ in which case $\cup F_\alpha$ is
a maximal good subset of $S$. Consider the set $S=\{(0,0,0), (0,0,1),
(1,1,0), (1,1,1)\}$ Here $R_1 = F_1=\{(0,0,0), (0,0,1)\}$ $R_2 =
F_2=\{(1,1,0), (1,1,1)\}$. But $F_1\cup F_2$ is not good. 

\item If bounded functions $f$ on $M$ have bounded solution of 
\begin{equation*}
f = u_1 + u_2 + \cdots + u_n, \tag{$*$}
\end{equation*}
then clearly bounded functions $f$ in $U(S)$ will have bounded solution
of ($*$).)

\item Suppose $X_1, X_2, \ldots , X_n$ are compact topological spaces.
Let $S\subset X_1\times\ldots \times X_n$ be compact and $M\subset S$ a
maximal good set. If continuous functions $f$ on $M$ have continuous
solution of ($*$) then clearly continuous function in $U(S)$ have
continuous solution of ($*$). 

\item Let $M\subset S$ is maximal good set and let $D$ be a boundary of
$M$ (refer \S~4 of \cite{KNB2}, for the definition of boundary).
Let $f$ be in $U(S)$. Let $g_i$ be defined on $D\cap \Pi_iS= D\cap
\Pi_iM$ for $1\leq i\leq n$. Then there exists a unique set
$v_1,\ldots,v_n$ such that $f=v_1+\cdots +v_n$ on $S$ and $v_i|_{D\cap
\Pi_iS}=g_i$. 

\item For any choice of $F_\alpha$'s ($F_\alpha$ as above), $\cup
F_\alpha$ is a maximal good set of $S$ if and only if $L$ is a loop in
$S$ implies $L\subset R_\alpha$ for some $\alpha$.\vspace{-2pc}
\end{enumerate}
\end{proof}

\begin{proof}
Suppose $\{x_1,\ldots,x_k\}$ is a loop in $S$, and
$\{x_1,\ldots,x_k\}\cap R_\alpha\neq \emptyset$ for more than one $\alpha$.
Then each such intersect is good because it is part of a loop. Let
$F_\alpha$ be a maximal full set in $R_\alpha$ such that $
\{x_1,\ldots,x_k\}\cap R_\alpha \subset F_\alpha$. (A good subset of a
set is contained in a maximal good subset of that set.) Then $\cup
F_\alpha$ contains a loop so is not good. Conversely suppose that any
given loop in $S$ is contained in some $R_\alpha$. Then $\cup F_\alpha$
cannot contain a loop because if it contains a loop then it is in some
$R_\alpha$, hence in some $F_\alpha$. But $F_\alpha$ cannot contain a\break
loop.
\end{proof}

\section*{Acknowledgements}

The author would like to thank Prof.~M~G~Nadkarni for suggesting the
problems and for encouragement and fruitful discussions. The author also
thanks Institute of Mathematical Sciences, for short visiting
appointments which made this work possible.

\end{document}